\documentclass[12pt,a4paper,reqno,twoside]{amsart}

\usepackage{amssymb}
\usepackage{amsmath}
\usepackage{amsthm}
\usepackage{amsfonts}
\usepackage{graphicx}
\usepackage{psfrag,color}
\usepackage[colorlinks,citecolor=blue,urlcolor=blue]{hyperref}
\usepackage{cite}
\usepackage[hmargin=2.5cm,vmargin=2.5cm]{geometry}

\newtheorem{theorem}{Theorem}[section]

\newtheorem{cor}{Corollary}[section]
\newtheorem{lem}{Lemma}[section]

\newtheorem{remark}{Remark}[section]

\newcommand{\beq}{\begin{equation}}
\newcommand{\eeq}{\end{equation}}
\newcommand{\beqn}{\begin{equation*}}
\newcommand{\eeqn}{\end{equation*}}

\newcommand{\slot}{\,\cdot\,}

\newcommand\dm{\,\rd\fm}

\newcommand\cB{{\mathcal B}}

\newcommand\cH{{\mathcal H}}

\newcommand\cL{{\mathcal L}}

\newcommand\bN{{\mathbb N}}

\newcommand\bP{{\mathbb P}}

\newcommand\bR{{\mathbb R}}
\newcommand\bS{{\mathbb S}}

\newcommand\rd{{\mathrm d}}

\newcommand\fm{{\mathfrak m}}

\def\bfone{\mathbf{1}}

\newcommand{\Ti}{(T^n)^{-1}_i}

\title[Explicit correlation bounds for expanding maps using coupling]{Explicit correlation bounds for expanding circle maps using the coupling method}
\author[Henri Sulku]{Henri Sulku}
\address[Henri Sulku]{
Department of Mathematics and Statistics, P.O.\ Box 68, Fin-00014 University of Helsinki, Finland.}
\email{henri.sulku@helsinki.fi}


\begin{document}

\begin{abstract}
In this paper, several fundamental facts, especially the existence and uniqueness of an absolutely continuous ergodic measure with an exponential mixing rate, are derived for smooth expanding circle maps. Although the results are classical, the coupling method is relatively modern and employed here for the first time to yield explicit bounds in terms of system constants. The work constitutes a part of the author's Bachelor's thesis. The author hopes that the manuscript can serve as a useful introduction to the flexible coupling method in the theory of dynamical systems.
\end{abstract}

\maketitle

\subsection*{Acknowledgements}

The thesis was written under a supervision of \mbox{dr.~Mikko~Stenlund} and I use this opportunity to express my gratitude for his support and professional guidance during the process.

\section{Introduction}

We assume that $T:\bS^1\to\bS^1$ is a smooth uniformly expanding map on the circle. More precisely, let $\bS^1$ denote the interval $[0,1]$ with its endpoints identified, let $T$ be $C^2$, and assume there exists $\lambda>1$ such that
\beq\label{eq:expansion0}
|T'(x)| \geq \lambda \qquad \forall\,x\in \bS^1.
\eeq
Our goal is to prove that there exists an invariant measure\footnote{A natural choice for the governing $\sigma$-algebra is the Borel $\sigma$-algebra $\cB$.} which is

\medskip
\noindent (i) absolutely continuous with respect to the Lebesgue measure and
\smallskip
\\
\noindent (ii) mixing (thus also ergodic).
\medskip

\noindent Recall that a measure is invariant if $\int_X f\circ T\,\rd\mu = \int_X f\,\rd\mu$ for all integrable functions~$f$. It is called mixing, provided that
\beqn
\int_X f\circ T^n g\,\rd\mu \longrightarrow \int_Xf\,\rd\mu\int_Xg\,\rd\mu\qquad \forall~ f,g\in L^2(\mu).
\eeqn
Note that here $\int_X f\circ T^n\,\rd\mu = \int_X f\,\rd\mu$ due to invariance. Thus, the mixing condition above has the probabilistic interpretation that the random variables $f\circ T^n(x)$ and $g(x)$ become asymptotically de-correlated as $n\to\infty$, if $x$ is distributed according to $\mu$.

\medskip
\noindent Actually, we want to prove stronger statements than these and obtain quantitative information about the invariant measure. In addition to (i) and (ii), we are going to show that the invariant density (the density of the invariant measure) is actually Lipschitz continuous and that it is exponentially mixing in a suitable class of observables; namely, the class of H\"older continuous functions of arbitrary order. Our main result is thus the following:

\begin{theorem}\label{thm_mainthm}
An invariant measure $\mu$ satisfying properties (i) and (ii) above exists. In fact, $\mu$ is the only ergodic measure satisfying (i). Additionally, for bounded $f$ and H\"older continuous $g$,
\beq
\bigg|\int_X f\circ T^n g\,\rd\mu - \int_Xf\,\rd\mu\int_Xg\,\rd\mu\bigg| \leq C\|f\|_\infty(\|g\|_\infty + H_\alpha(g))\theta^{\alpha n},
\eeq
with
\beq
C = 96(2 + \Omega)^2
\eeq
and
\beq
\theta = (1 - \exp(-3(\Omega+1))^{(\log(\lambda)/(4(\Omega+1)))} \in (0,1).
\eeq
\end{theorem}

\noindent Here $H_\alpha(g)$ stands for the H\"older constant of a function $g$.

\medskip
\noindent The treatise, and the proof of \ref{thm_mainthm} as well, is divided into five sections: first we express some preliminary results which will play a central role in the study of the system. Then, in the third section, we prove the existence of an invariant measure with a Lipschitz continuous density. The coupling argument is employed in the fourth section where we prove the exponential decay of correlation coefficients for H\"older continuous observables. Some of the properties are studied in a bit more general setting in \cite{OSY} as well. Finally, at the fifth section we prove the mixing property by extending the correlation decay from $C^\alpha$ to $L^1(\mu)$. Remark that the exponential rate is not preserved but even the continuous functions can be too irregular instead. Additionally, a couple of fundamental results are expressed in the appendix because those are somewhat isolated from the general storyline.

\section{Preliminaries}\label{sec:prelim}
We begin by introducing some basic facts that will be used throughout the text. The expansion property below is a manifestation of the hyperbolic nature of the dynamics and inverse branches are needed to treat the transfer operator effectively. The transfer operator itself will be a crucial tool in the analysis of the dynamics.

\subsection{Expansion properties and inverse branches}
First, since $T$ is $C^2$ and \eqref{eq:expansion0} holds, $T'$ cannot change sign. Therefore, we can and do assume that~$T'$ is strictly positive:
\beq\label{eq:expansion}
T'(x) \geq \lambda \qquad \forall\,x\in \bS^1.
\eeq

\noindent We denote the standard metric on $\bS^1$ by $d(\slot,\slot)$, i.e., $d(x,y) = \min\{|x-y|,1-|x-y|\}$ and arc lengths by~$|\slot|$. Observe that, given an arc $I\subsetneq \bS^1$ such that $T(I)\neq \bS^1$, it is obviously true by the expanding property~\eqref{eq:expansion} that  $T(I)$ is an arc with $|T(I)|\geq \lambda | I |$.
Notice that the quantity $w=\int_{[0,1]} T'(x)\,\rd x$ is a positive integer. We present the following, rather obvious, result without proof.

\begin{lem}\label{lem:branches}
Let $J\subsetneq\bS^1$ be an arc. Then the map $T$ has exactly $w$ well-defined inverse branches on $J$. More precisely, there exist arcs $I_i\subsetneq\bS^1$ and maps $T_i^{-1}:J\to I_i$, $1\leq i\leq w$, such that each restriction $T | I_i$ is a one-to-one and onto $J$ (in fact a $C^2$ diffeomorphism in the interior of $I_i$) and $T_i^{-1}$ is its inverse.
\end{lem}
\noindent By the earlier observation, $|J| = |T(I_i)|\geq \lambda |I_i|$. Therefore, Lemma~\ref{lem:branches} can be applied repeatedly, resulting in exactly $w^n$ well-defined inverse branches of $T^n$ on $J$. Let us denote them $(T^n)_i^{-1}:J\to I_i^n$, $1\leq i\leq w^n$, where each $I_i^n$ is an arc. Iterating the previous bound, $|I_i^n| \leq \lambda^{-n}|J|$. This results in the lemma below.
\begin{lem}\label{lem:backward_contraction}
Let $x,y\in\bS^1$ be arbitrary. Suppose $J\subset \bS^1$ is an arc with $|J|\leq \tfrac12$ and $x,y\in J$. Then
\beqn
d\bigl((T^n)_i^{-1}(x),(T^n)_i^{-1}(y)\bigr) \leq \lambda^{-n} \,d(x,y)
\eeqn
for all $1\leq i\leq w^n$ and all $n\geq 1$. (Here $(T^n)_i^{-1}$ are the inverse branches of $T^n$ on $J$.)
\end{lem}
\begin{proof}
Consider the arc $J'\subset J$ with endpoints $x,y$. It is a subset of a semicircle in~$\bS^1$, so that $d(x,y) = |J'|$. Hence, $d\bigl((T^n)_i^{-1}(x),(T^n)_i^{-1}(y)\bigr) = |(T^n)^{-1}_i(J')|\leq \lambda^{-n} | J' | = \lambda^{-n} d(x,y)$ by the observation preceding the lemma.
\end{proof}

\subsection{Transfer operator}
The transfer operator $\cL_T:L^1(\fm)\to L^1(\fm)$ is defined by the ``duality equation'':
\beq
\label{eq:transfer_dual}
\int_X f\circ T\cdot g\dm = \int_X f\cdot \cL_Tg\dm \qquad\forall~f\in L^\infty(\fm), g\in L^1(\fm).
\eeq
In particular, if $\psi$ is the density of a probability measure $\nu$, $\cL\psi$ is the density of the push-forward measure $T_*\nu = \nu\circ T^{-1}$. In what follows, the $T$-subscript is omitted if there is no risk of misconception and $\cL$ denotes the transfer operator associated with $T$. One can easily confirm that $\cL u$ has the expression
\beq
\label{eq:clrep}
\cL u(x) = \sum_{y\in T^{-1}\{x\}} \frac{u(y)}{T'(y)}.
\eeq
The missing details in the proof of the following lemma, which summarises the basic properties of $\cL$, are left to the reader.

\begin{lem}
\label{lem:Trans_prop}
For transfer operator related to expanding mapping $T$,
\begin{enumerate}
\item $|\cL u| \leq \cL|u|$ for all $u:\bS^1\to\bR$,
\item $0 \leq u \leq v \Longrightarrow 0 \leq \cL u\leq \cL v$ for all $u,v:\bS^1\to\bR$,
\item $\int_X \cL u\dm = \int_X u\dm$ for all $u\in L^1(\fm)$~~\textrm{and}
\item $\|\cL u\|_{L^1(\fm)} \leq \|u\|_{L^1(\fm)}$ for all $u\in L^1(\fm)$.
\end{enumerate}
\end{lem}

\smallskip
\begin{proof}~\\
1. Use \eqref{eq:clrep} and triangle inequality.\\
2. Use \eqref{eq:clrep} and linearity of $\cL$.\\
3. By Equation \eqref{eq:transfer_dual}
\beqn
\int_X\cL u\,\rd\fm = \int_X1\cdot\cL u\,\rd\fm = \int_X1\circ T\cdot u\,\rd\fm = \int_X u\,\rd\fm~~\forall~u\in L^2
\eeqn
4. Combine properties 1 and 3.
\end{proof}

\noindent One can also deduce expressions
\begin{align}
\cL u(x) &= \sum^{w}_{y=1} \frac{u(T^{-1}_{i}(x))}{T'(T^{-1}_{i}(x))} \\
(\cL u)'(x) &= \sum^{w}_{i=1}\left(\frac{u'(T^{-1}_{i}(x))}{T'(T^{-1}_{i}(x))^2}-\frac{u(T^{-1}_{i}(x))}{T'(T^{-1}_{i}(x))^3}T''(T^{-1}_{i}(x))\right) \label{eq:Dtransfer}
\end{align}
in case $u\in C^1$ using inverse branches of $T$.

\section{Existence of a.c.i.m.}

The existence of the a.c.i.m. is proved by a compactness argument similar to Arzel\`a--Ascoli theorem. More precisely, we are about to prove that $C^1$-norms of the pushforward densities remain uniformly bounded and bounded subsets of $C^1$ are compact in the space of Lipschitz functions endowed with its natural norm. It is worth mentioning that the strategy we use, that is, establishing invariance by taking averages along trajectories and using compactness to deduce that a limit exists, is widely used method to construct invariant measures (see e.g.~the proof of  Krylov--Bogolyubov Theorem \cite[p.~151]{Wal00}).

\subsection{Boundedness}
In this section, we derive estimates for the transfer operator which C. Liverani motivated in \cite{Liv}. The ultimate result we wish to prove is that for every fixed $\psi\in C^1$ there is $L>0$ such that
\beq
\label{eq:thegoal1}
\sup_{n\in\bN}\|\cL^n\psi\|_{C^1} = \sup_{n\in\bN}\left(\|\cL^n\psi\|_\infty+\|(\cL^n\psi)'\|_\infty\right) \leq L.
\eeq
We begin by pointing out that the frequently used estimate
\beq
\label{eq:Destimate}
|(\cL\psi)'| \leq \frac{1}{\lambda}\cL|\psi'|+\frac{\|T''\|_\infty}{\lambda^2}\cL|\psi|
\eeq
follows directly from \eqref{eq:Dtransfer}. Using this estimate and the fact the $\cL$ is a contraction in $L^1(\fm)$ (property 4 in Lemma \ref{lem:Trans_prop}) we can control $L^1$-norms of the pushforwards. Let us prove that these are uniformly bounded with respect to $n$.
\begin{lem}
Let $\psi\in C^1$ be arbitrary. For every $n\in\bN$ we have the following upper bound
\beq
\label{eq:unif_L1bound}
\|(\cL^n\psi)'\|_{L^1(\fm)} \leq \frac{1}{\lambda^n}\|\psi'\|_{L^1(\fm)}+\frac{\|T''\|_\infty}{\lambda(\lambda-1)}\|\psi\|_{L^1(\fm)}.
\eeq 
\end{lem}

\begin{proof} The $n=1$ case is a direct consequence of the contracting nature of $\cL$ and \eqref{eq:Destimate}:
\beqn
\|(\cL\psi)'\|_{L^1(\fm)} \leq \frac{1}{\lambda}\|\cL|\psi'|\|_{L^1(\fm)}+\frac{\|T''\|_\infty}{\lambda^2}\|\cL|\psi|\|_{L^1(\fm)} \leq \frac{1}{\lambda}\|\psi'\|_{L^1(\fm)}+\frac{\|T''\|_\infty}{\lambda(\lambda-1)}\|\psi\|_{L^1(\fm)}.
\eeqn
Next we claim that
\beq
\label{eq:Lnuestimate}
\|(\cL^n \psi)'\|_{L^1(\fm)} \leq \frac{1}{\lambda^n}\|\psi'\|_{L^1(\fm)}+\|T''\|_{\infty}\|\psi\|_{L^1(\fm)}\sum^{n+1}_{i=2}\frac{1}{\lambda^i}
\eeq
holds for any $n \in \bN$. To this end, we proceed by induction. Suppose \eqref{eq:Lnuestimate} holds for $n=k$. For $n=k+1$ we have
\begin{align}
\|(\cL^{k+1}\psi)'\|_{L^1(\fm)} &\leq \frac{1}{\lambda}\|\cL|(\cL^k \psi)'|\|_{L^1(\fm)}+\frac{\|T''\|_\infty}{\lambda^2}\|\cL|\cL^{k}\psi|\|_{L^1(\fm)} \nonumber \\
&\leq\frac{1}{\lambda}\|(\cL^k\psi)'\|_{L^1(\fm)}+\frac{\|T''\|_\infty}{\lambda^2}\|\psi\|_{L^1(\fm)} \nonumber \\
&\leq \frac{1}{\lambda}\left(\frac{1}{\lambda^k}\|\psi'\|_{L^1(\fm)}+\|T''\|_{\infty}\|\psi\|_{L^1(\fm)}\sum^{k+1}_{i=2}\frac{1}{\lambda^i}\right)+\frac{\|T''\|_\infty}{\lambda^2}\|\psi\|_{L^1(\fm)} \nonumber \\
&\leq\frac{1}{\lambda^{k+1}}\|\psi'\|_{L^1(\fm)}+\|T''\|_{\infty}\|\psi\|_{L^1(\fm)}\sum^{k+2}_{i=2}\frac{1}{\lambda^i} \nonumber.
\end{align}
Thus, we conclude that \eqref{eq:Lnuestimate} holds for any $n \in \bN$. Finally, we observe that
\beqn
\sum_{i=2}^{n+1}\frac{1}{\lambda^i} \leq \frac{1}{\lambda(\lambda-1)}~~\forall~n\in\bN
\eeqn
which finishes the proof.
\end{proof}

\medskip
\noindent
The constant $\|T''\|_\infty/(\lambda(\lambda-1))$ appearing in the left-hand side of \eqref{eq:unif_L1bound} is so extensively used throughout this study that we denote it by the symbol $\Omega$, i.e.,
\beq
\Omega := \frac{\|T''\|_\infty}{\lambda(\lambda-1)}.
\eeq
Constant $\Omega$ in a way measures the ''regularity'' of the dynamics, that is, the linearity of the dynamics compared to the minimum magnitude of stretching.

\begin{remark}
The contribution of the derivative of the initial funtion $\psi$ to the upper bound in \eqref{eq:unif_L1bound} becomes neglible after a sufficiently long time. This is a very general property of expanding dynamics; in the long term, sufficiently regular functions become even more smoothly distributed.
\end{remark}

\noindent The reason to examine $L^1$-norms is that in case of $C^1$ probability densities we can bound sup-norms by $L^1$-norms. More precisely, if $\psi$ is differentiable density there is at least one point $x_0$ in which $\psi(x_0)=1$. Otherwise it would not be possible to have $\int_X\psi\dm = 1$. Thus,
\beqn
|\psi(x)| = |\int^{x}_{x_0}\psi'(t)dt+\psi(x_0)| \leq \int^{x}_{x_0}|\psi'(t)|dt+1 \leq \|\psi'\|_{L^1(\fm)}+1
\eeqn
holds for every $x\in\bS^1$ and therefore
\beq
\label{eq:L1estimate}
\|\psi\|_\infty \leq \|\psi'\|_{L^1(\fm)}+1.
\eeq
Remark that $\|\cL^n v\|_\infty \leq \|v\|_\infty\|\cL^n\bfone\|_\infty$ holds by the second part of lemma \ref{lem:Trans_prop}, where $\bfone = \chi_{\bS^1} \in C^1$. We now have the tools to establish the desired boundedness result.

\begin{theorem}
\label{thm:sup_bound}
For every bounded $\psi$
\beqn
\sup_{n\in\bN}\|\cL^n\psi\|_\infty \leq (1+\Omega)\|\psi\|_\infty
\eeqn
\end{theorem}

\begin{theorem}
\label{thm:supD_bound}
For every continuously differentiable $\psi$
\beqn
\sup_{n\in\bN}\|\cL^n\psi\|_{C^1} \leq (1+\Omega)^2\|\psi\|_{C^1}
\eeqn
\end{theorem}

\begin{proof}[Proof of the Theorem \ref{thm:sup_bound}]
Notice that $\cL^n\bfone$ is $C^1$ for all $n\in\bN$. Therefore, estimates \eqref{eq:unif_L1bound} and \eqref{eq:L1estimate} results in
\beqn
\|\cL^n\psi\|_\infty \leq \|\psi\|_\infty\|\cL^n\bfone\|_\infty \leq \|\psi\|_\infty\left(\|(\cL^n\bfone)'\|_{L^1(\fm)}+1\right) \leq \|\psi\|_\infty\left(\frac{\|T''\|_\infty}{\lambda(\lambda-1)}+1\right).
\eeqn
\end{proof}

\begin{proof}[Proof of the Theorem \ref{thm:supD_bound}]
In order to achieve uniform bound, we first claim that
\beqn
|(\cL^n \psi)'| \leq \frac{1}{\lambda^n}\cL^n|\psi'|+\|T''\|_\infty\cL^{n}|\psi|\sum^{n+1}_{i=2}\frac{1}{\lambda^i}.
\eeqn
By estimate \eqref{eq:Destimate} this holds when $n=1$. For $n>1$ we proceed by induction. Suppose the inequality holds for $n=k$. If $n=k+1$ it holds that
\begin{align}
|(\cL^{k+1}\psi)'| &\leq \frac{1}{\lambda}\cL|(\cL^k\psi)'|+\frac{\|T''\|_\infty}{\lambda^2}\cL^{k+1}|\psi| \nonumber
\\
&\leq \frac{1}{\lambda}\cL\bigg|\frac{1}{\lambda^k}\cL^k|\psi'|+\|T''\|_\infty\cL^{k}|\psi|\sum^{k+1}_{i=2}\frac{1}{\lambda^i}\bigg|+\frac{\|T''\|_\infty}{\lambda^2}\cL^{k+1}|\psi|\nonumber
\\
&\leq \frac{1}{\lambda^{k+1}}\cL^{k+1}|\psi'|+\|T''\|_\infty\cL^{k+1}|\psi|\sum^{k+1}_{i=2}\frac{1}{\lambda^{i+1}} +\frac{\|T''\|_\infty}{\lambda^2}\cL^{k+1}|\psi|\nonumber
\\
&=\frac{1}{\lambda^{k+1}}\cL^{k+1}|\psi'|+\|T''\|_\infty\cL^{k+1}|\psi|\sum^{k+2}_{i=2}\frac{1}{\lambda^i}. \nonumber
\end{align}
Therefore, by the induction principle, the estimate holds for any $n\in\bN$. From this estimate it is easy to deduce the following bound
\beq
\|(\cL^n \psi)'\|_{\infty} \leq (1+\Omega)(\|\psi'\|_\infty+\Omega\|\psi\|_\infty)
\eeq
\noindent What we have proved is that
\begin{align}
\label{eq:uniform_bound}
\|\cL^n\psi\|_{C^1} &\leq (1+\Omega)^2\|\psi\|_\infty+(1+\Omega)\|\psi'\|_\infty\nonumber\\
&\leq (1+\Omega)^2\|\psi\|_\infty+(1+\Omega)^2\|\psi'\|_\infty\nonumber \\
&= (1+\Omega)^2\|\psi\|_{C^1}
\end{align}
for all natural numbers $n$ and continuously differentiable $\psi$.
\end{proof}

\subsection{Compactness}

We begin by proving the compactness result that guarantees the existence of a Lipschitz continuous invariant density.

\begin{lem}
Suppose that a sequence of $C^1$ functions $f_n:\bS^1\to\bR$ is bounded, i.e., $\sup_n\|f_n\|_{C^1}\leq L<\infty$. Then it has a subsequence $(f_{n_k})_k$ which converges uniformly to a Lipschitz function $f$ having Lipschitz constant $\mathrm{Lip}(f) \leq L$. 
\end{lem}

\begin{proof}
The existence of the uniformly convergent subsequence is just an application of the Arzel\`a--Ascoli theorem: $C^1$--bounded sequence is, of course, equicontinuous and pointwise uniformly bounded so it satisfies the conditions of the theorem. Since the uniformly convergent exists, the proof is finished by noting that uniform limit of functions with uniformly bounded Lipschitz constants is Lipschitz continuous.
\end{proof}

\noindent The lemma can be generalised a little and re-interpreted as follows: Closed, bounded, subsets of $C^1$ are compact in the space of Lipschitz continuous functions endowed with its natural norm $\|\cdot\| = \|\cdot\|_\infty+\mathrm{Lip}(\cdot)$. With the aid of this lemma we are ready to prove the main theorem of the section.

\smallskip
\begin{theorem}[The existence of the Lipschitz continuous invariant density]~
\label{existence_density}
\noindent There is an invariant measure $\mu$ such that it is absolutely continuous w.r.t.~Lebesgue measure and its density is Lipschitz continuous with its Lipschitz constant bounded by a quantity depending only on the dynamics of the system. More precisely, the density $\phi = \rd\mu/\rd\fm$ satisfies $\mathrm{Lip}(\phi) \leq (1+\Omega)^2$.
\end{theorem} 

\begin{proof}
Let $\psi$ be $C^1$-density associated to absolutely continuous, not necessarily invariant, probability measure $\nu$. By theorem \ref{thm:supD_bound}, we know that densities $\cL^n\psi$ are contained in some closed ball in $C^1$. Now, define a sequence of probability measures $\mu_n$ by densities of the form
\beq\label{eq:orbit_average}
\psi_N = \frac{1}{N}\sum^{N-1}_{n=0}\cL^n\psi.
\eeq
Since the averages are contained in the closed ball, we know that there is a subsequence $(\psi)_{n_j}$ such that $\psi_{n_j}$ converges uniformly to a Lipschitz density $\phi$. Let $\mu$ be the probability measure defined by the limit density $\phi$. In order to show it is invariant, let $A$ be an arbitrary measurable set. By the bounded convergence theorem and trivial observation that
\beqn
\int_{T^{-1}A} \psi_{n_j}\,\rd\fm = \int_{\bS^1}1_{T^{-1}A}\, \psi_{n_j}\,\rd\fm = \int_{\bS^1}1_{A}\circ T\, \psi_{n_j}\,\rd\fm =  \int_{\bS^1}1_{A}\, \cL\psi_{n_j}\,\rd\fm = \int_{A} \cL\psi_{n_j}\,\rd\fm,
\eeqn
we have the following equality
\beqn
\begin{split}
\mu(T^{-1}A) - \mu(A) 
& =  \int_{T^{-1}A} \lim_{j\to\infty}\psi_{n_j} \,\rd\fm - \int_A \lim_{j\to\infty}\psi_{n_j} \,\rd\fm
\\
& =  \lim_{j\to\infty}\int_{T^{-1}A}\psi_{n_j} \,\rd\fm -  \lim_{j\to\infty}\int_A \psi_{n_j} \,\rd\fm
\\
& = \lim_{j\to\infty}\int_A (\cL\psi_{n_j}-\psi_{n_j})\,\rd\fm.
\end{split}
\eeqn
Now the proof of invariance is accomplished by the following estimate:
\begin{align}
\left|\int_A (\cL\psi_{n_j}-\psi_{n_j}) \, \rd\fm\right| & = \left |\int_A\frac{\cL^{n_j}\psi}{n_j}-\frac{\psi}{n_j}\rd\fm \right|\nonumber \\
&\leq \frac{\urladdr{}1}{n_j}(\|\cL^{n_j}\psi\|_{L^1}+\|\psi\|_{L^1}) \nonumber \\
&\leq \frac{2\|\psi\|_{L^1}}{n_j} \longrightarrow 0 .\nonumber
\end{align}

\medskip
\noindent An upper bound for Lipschitz constant is obtained by estimates on the Lipschitz constants of $\psi_{n_j}$:
\begin{align}
\mathrm{Lip}(\phi) &\leq \sup_{j\in\bN}\mathrm{Lip}(\psi_{n_j}) \nonumber \\
&\leq \sup_{j\in\bN}\|\psi_{n_j}\|_{C^1} \nonumber \\
&\leq (1+\Omega)^2\|\psi\|_{C^1} \nonumber
\end{align}
Additionally, since $\psi$ was arbitrary, we may pick $\psi \equiv 1$ to obtain
\beq
\mathrm{Lip}(\phi) \leq (1+\Omega)^2.
\eeq
\end{proof}
\noindent Remark that this is just an existence theorem, not uniqueness; we have showed that \emph{every}~$C^1$ density $\psi$ results in \emph{at least one} invariant density. Nevertheless, uniqueness follows from the convergence of an arbitrary density towards the invariant density.

\section{Exponential decay of correlation coefficients}
\label{consid_page}
Recall that the mixing property can be equivalently stated in terms of decay of \emph{correlation coefficients} of $L^2$-functions:
\beq\label{eq:mixing}
\lim_{n\to\infty}\int f\circ T^n\cdot g\,\rd\mu = \int f\,\rd\mu \int g\,\rd\mu \qquad\forall \,f,g\in L^2(\mu).
\eeq
In previous section, we proved that $\mu$ is absolutely continuous with respect to the Lebesgue measure. Denoting the density by $\phi$, mixing is equivalent to
\beqn
\lim_{n\to\infty}\int f\circ T^n\cdot g\phi\,\rd\fm = \int f\phi\,\rd\fm \int g\phi\,\rd\fm \qquad\forall \,f,g\in L^2(\mu).
\eeqn

\medskip
\noindent As motivated at the beginning of the paper, we are willing to restrict ourselves to H\"older continuous $g$ in order to obtain concrete estimates on the error terms
\beqn
\bigg|\int f\circ T^n\cdot g\phi\,\rd\fm - \int f\phi\,\rd\fm \int g\phi\,\rd\fm\bigg|.
\eeqn
To this end, remark that by the linearity of integrals, we may add $2\int f\phi\,\rd\fm\|g\|_\infty$ both sides and normalise the equality to
\beqn
\bigg|\int g\,\rd\mu+2\|g\|_\infty\bigg| \cdot \bigg|\int f\circ T^n\cdot \frac{\phi(g+2\|g\|_\infty)}{\int g\,\rd\mu+2\|g\|_\infty}\,\rd\fm - \int f\phi\,\rd\fm\bigg|.
\eeqn
Denoting
\beqn
\psi =  \frac{\phi(g+2\|g\|_\infty)}{\int g\,\rd\mu+2\|g\|_\infty},
\eeqn
we can reduce the decay of correlation coefficients to $L^1(\fm)$--convergence of probability densities towards the invariant density:
\beqn
\bigg|\int f\circ T^n\cdot g\phi\,\rd\fm - \int f\phi\,\rd\fm \int g\phi\,\rd\fm\bigg| \leq 3\|g\|_\infty\|f\|_\infty\|\cL^n\psi-\phi\|_{L^1(\fm)}.
\eeqn

\medskip
\noindent Therefore, it is sufficient to prove the inequality
\beq\label{exponential_estimate}
\|\cL^n\psi-\phi\|_{L^1(\fm)} \leq D\theta^{n}\qquad\forall~\psi\in C^\alpha
\eeq
with $D > 0$, $\theta \in (0,1)$, in order to achieve exponentially decaying estimates. This will be proved after we have developed some tools needed to construct the coupling argument.

\medskip
\noindent In what follows, we use the notation
\beqn
H_\alpha(f) = \sup_{x,y\in X}\frac{|f(x)-f(y)|}{d(x,y)^\alpha}
\eeqn
for the H\"older coefficient of a function $f$ and denote the class of H\"older continuous functions of order $\alpha$ by $C^\alpha$. That is,
\beqn
C^\alpha = \{f:X\to\bR: H_\alpha(f)<\infty\}
\eeqn

\subsection{Preliminary results}

The following \emph{distortion bound} is central for understanding the structure of push-forward densities $\cL^n\psi$. The core idea of the proof is basically the same as Lai-Sang Young's in \cite[p.~33-34]{LSY}.

\begin{lem}\label{lem:distortion}
Let $n\in\bN$ be arbitrary. For any $x,y\in\bS^1$,
\beq
e^{-\Omega d(x,y)}\leq \frac{(T^n)'((T^n)^{-1}_ix)}{(T^n)'((T^n)^{-1}_iy)} \leq e^{\Omega d(x,y)},
\eeq
in which $\Omega = \|T''\|_\infty/(\lambda(\lambda-1))$ \emph{is independent of $n$}. Here $\Ti$ is the $i$th branch of the inverse of $T^n$ on a given arc $J\subset \bS^1$ of length $|J|\leq \frac12$ containing both $x$ and $y$.
\end{lem}

\begin{proof}
Let $J$ be as in the statement of the lemma. Lemma~\ref{lem:backward_contraction} implies that,  for an arbitrary $k\geq 1$, the inverse branches of $T^k$ on $J$ are well defined. For brevity, let $x_{-k}$ and $y_{-k}$ denote the preimages of $x$ and $y$, respectively, along the same branch.

\medskip
\noindent Since the logarithm is $\lambda^{-1}$--Lipschitz in $[\lambda,\|T'\|_\infty]$, we can estimate
\begin{align*}
\log \frac{(T^n)'(x_{-n})}{(T^n)'(y_{-n})} &\leq |\log((T^n)'(x_{-n}))-\log((T^n)'(y_{-n}))| \\
&\leq \sum^{n-1}_{i=1}|\log(T'(T^ix_{-n}))-\log(T'(T^iy_{-n}))| \\
&\leq \sum^{n-1}_{i=1}\frac{1}{\lambda}\|T''\|_\infty d(x_{-n+i},y_{-n+i}) \\
&\leq \sum^{n-1}_{i=1}\frac{1}{\lambda}\|T''\|_\infty\lambda^{-n+i}d(x,y) \\
&\leq \frac{\|T''\|_\infty}{\lambda(\lambda-1)}d(x,y) = \Omega d(x,y).
\end{align*}
A similar estimate is obtained by interchanging $x$ and $y$, which proves the claim.
\end{proof}

\medskip
\noindent This lemma is immediately used in the following theorem which is the cornerstone of our coupling argument. The trick is that sometimes, particularly in our case, it is more convenient to work with logarithms of functions instead of functions themselves.

\begin{theorem}
\label{thm:Hoelder_estimate}
Suppose $\psi$ is a strictly positive probability density and that $\log\psi \in C^\alpha$.
Then $\cL^n\psi$ has the same properties for every $n\in\bN$ and 
\beqn
H_\alpha(\log\cL^n\psi) \leq \frac{H_\alpha(\log\psi)}{\lambda^{\alpha n}}+ \Omega.
\eeqn
\end{theorem}

\begin{proof}
Let $J\subset \bS^1$ be an interval with $|J|\leq \frac12$.
Given an initial probability density $\psi$, we introduce the notation
\beqn
\psi_{n,i}(x) = \frac{\psi((T^n)^{-1}_i x)}{(T^n)'(\Ti x)} \ ,\qquad x\in J,
\eeqn
where $(T^n)^{-1}_i$ denotes the $i$th inverse branch of $T$ on $J$. Then
\beqn
\cL^n\psi(x) = \sum_{i=1}^{w^n} \psi_{n,i}(x).
\eeqn

The number $w$ has been introduced in Section~\ref{sec:prelim}.

\medskip
\noindent Next, let $x,y\in\bS^1$ be arbitrary. Without loss of generality, we may assume both points belong to $J$. Therefore,
\beqn
\begin{split}
\left|\log\frac{\psi_{n,i}(x)}{\psi_{n,i}(y)}\right| 
&
\leq \left|\log\frac{\psi((T^n)^{-1}_i x)}{\psi((T^n)^{-1}_i y)}\right| + \left|\log\frac{(T^n)'(\Ti y )}{(T^n)'(\Ti x)}\right|
\\
&
\leq H_\alpha(\log\psi) \,d((T^n)^{-1}_i x, (T^n)^{-1}_i y)^\alpha+ \Omega d(x,y)
\\
&
\leq \lambda^{-\alpha n} H_\alpha(\log\psi) \,d(x,y)^\alpha+ \Omega d(x,y)^\alpha.
\end{split}
\eeqn
For brevity, denote the right side of the last bound by $R$. Then
\beqn
e^{-R} \leq \frac{\psi_{n,i}(x)}{\psi_{n,i}(y)} \leq e^{R},\textrm\qquad{ i.e. }\qquad e^{-R}\psi_{n,i}(y) \leq \psi_{n,i}(x)\leq e^{R}\psi_{n,i}(y).
\eeqn
Summing over $i$, we get
\beqn
e^{-R}\cL^n\psi(y) \leq \cL^n\psi(x)\leq e^{R}\cL^n\psi(y),\textrm\qquad{ i.e. }\qquad e^{-R} \leq \frac{\cL^n\psi(x)}{\cL^n\psi(y)} \leq e^{R}.
\eeqn
Taking logarithms, this yields
\beqn
\left|\log\frac{\cL^n\psi(x)}{\cL^n\psi(y)}\right| 
\leq (\lambda^\alpha)^{-n} H_\alpha(\log\psi) \,d(x,y)^\alpha+ \Omega d(x,y)^\alpha,
\eeqn
which is the desired bound.
\end{proof}

\medskip
\noindent We finish this section by demonstrating how estimates on H\"older constants of logarithms can be used to obtain information about approppriate \emph{probability} densities.

\medskip
\noindent As mentioned earlier, for every probability density $\psi$, there is a point $x_0$ such that $\psi(x_0) = 1$. Therefore,
\beqn
|\log(\psi(x))| = |\log(\psi(x))-\log(\psi(x_0))| \leq H_\alpha(\log(\psi))d(x,x_0)^\alpha \leq H_\alpha(\log(\psi)),
\eeqn
i.e.
\beq\label{eq:log_bounds}
\exp(-H_\alpha(\log(\psi))) \leq \psi(x) \leq \exp(H_\alpha(\log(\psi))),
\eeq
holds for every $x\in\bS^1$. Furthermore, for every $x,y \in \bS^1$
\beqn
d(x,y)^\alpha H_\alpha(\log\psi) \geq |\log(\psi(x))-\log(\psi(y))| = |\int^{\psi(x)}_{\psi(y)}\frac{dt}{t}| \geq \frac{1}{\|\psi\|_\infty}|\psi(x)-\psi(y)|
\eeqn
from which it follows that
\beq\label{eq:const_log_bound}
H_\alpha(\psi) \leq H_\alpha(\log\psi)\|\psi\|_\infty \leq H_\alpha(\log\psi)\exp\left(H_\alpha(\log\psi)\right).
\eeq

\noindent Therefore, H\"older continuity of the logarithm implies H\"oelder continuity of the density itself. Additionally, it results in upper and lower bounds on the density.

\subsection{Coupling}
We begin by proving the following theorem which summarises the relevant consequences of the results of the previous subsection. To this end, let us introduce the following notation:
\beqn
\cH_D = \{\psi:X\to\bR: \psi > 0, H_\alpha(\log\psi) \leq D, \|\psi\|_{L_1(\fm)} = 1\}.
\eeqn
$\alpha\in (0,1]$ is fixed and $D \in \bR_+$.

\medskip
\begin{theorem}\label{thm:positive_convergence}~

\smallskip
\begin{enumerate}
\item For all $B > 0$ there is $N = N(B)\in\bN$ such that for all $\psi\in\cH_B$ and $n > N$ $\cL^n\psi \in \cH_{\Omega+1}$
\item There is $a > 0$ such that $\psi \geq 2a > 0$ holds for every $\psi\in\cH_{\Omega+1}$.
\item There is $K$ such that $\tilde{\psi} := (\psi-a)/(1-a)\in\cH_K$ for all $\psi\in\cH_{\Omega+1}$. The constant $a$ is the one from the previous item.
\end{enumerate}
\end{theorem}

\begin{proof}~

(1) By theorem \ref{thm:Hoelder_estimate}, this holds when
\beqn
\frac{B}{\lambda^{\alpha N}}+\Omega < \Omega+1 \textrm{ i.e. } N > \frac{\log(B)}{\alpha\log(\lambda)}.
\eeqn
Thus, $N(B)$ can be chosen to be the integer part of $\log(B)/(\alpha\log(\lambda))+1$.

\medskip
(2) By equation \eqref{eq:log_bounds},
\beqn
\psi(x) \geq \exp\left(-H_\alpha(\log\psi)\right) \geq \exp\left(-(\Omega+1)\right) =: 2a
\eeqn
holds for every $\psi\in\cH_{\Omega+1}$.

\medskip
(3) By equations \eqref{eq:log_bounds} and \eqref{eq:const_log_bound},
\begin{align*}
&\quad |\log(\frac{\psi(x)-a}{1-a})-\log(\frac{\psi(y)-a}{1-a})| \\
&\leq |\psi(x)-\psi(y)|\sup_x\frac{1}{\psi(x)-a} \\
&\leq \frac{1}{a}H_\alpha(\log\psi)\exp\left(\frac{H_\alpha(\log\psi)}{2^\alpha}\right)d(x,y)^\alpha \\
&\leq 2(\Omega+1)\exp\left(2(\Omega+1)\right)d(x,y)^\alpha \\
&\leq \exp\left(4(\Omega+1)\right)d(x,y)^\alpha \\
&=: Kd(x,y)^\alpha
\end{align*}
holds for all $\psi\in\cH_{\Omega+1}$. Thus, $\tilde{\psi} \in\cH_K$ for all $\psi\in\cH_{\Omega+1}$. In particular,
\beqn
N(K) = \frac{4(\Omega+1)}{\alpha\log(\lambda)}.
\eeqn
\end{proof}
\noindent Now we are ready to proceed to the coupling argument itself. In what follows, let $\alpha$ and~$T$ be fixed and $K$ and $a$ be the constants from the previous theorem. Let us prove the exponential convergence first for functions in $\cH_K$.

\medskip
\noindent Let $\psi_1$ and $\psi_2$ be arbitrary densities from $\cH_K$. After time $N = N(K)$ calculated earlier the pushforwards of these densities lie in $\cH_{\Omega+1}$. Thus, for $n > N$
\beqn
\cL^n\psi_i = a + (1-a)\tilde{\psi}_i
\eeqn
for some unique $\tilde{\psi}_i\in\cH_K$. Using this equality, $\cL$ being linear contraction in $L^1(\fm)$ yields
\beqn
\|\cL^n\psi_1 - \cL^n\psi_2\|_{L^1(\fm)} \leq (1-a)\|\tilde{\psi}_1 - \tilde{\psi}_2\|_{L^1(\fm)}.
\eeqn
Now we can apply the same argument for $\tilde{\psi}_1$ and $\tilde{\psi}_2$ which belong to $\cH_K$ as well. Therefore, for every $n>2N$ the inequality
\beqn
\|\cL^n\psi_1 - \cL^n\psi_2\|_{L^1(\fm)} \leq (1-a)^2\|\tilde{\psi}^1_1 - \tilde{\psi}^1_2\|_{L^1(\fm)}.
\eeqn
holds for some $\tilde\psi^{1}_i\in\cH_K$. Moreover, by simple induction argument we deduce that
\beqn
n \in [kN,(k+1)N] \Longrightarrow \|\cL^n\psi_1 - \cL^n\psi_2\|_{L^1(\fm)} \leq (1-a)^k\|\tilde{\psi}^k_1 - \tilde{\psi}^k_2\|_{L^1(\fm)}
\eeqn
holds for every $k\in\bN$ with some unique $\tilde{\psi}^{k}_i\in\cH_K$. On the other hand, $n/N \leq (k+1)$ yields
\beqn
(1-a)^{k} \leq (1-a)^{n/N-1}.
\eeqn
Combining these results and uniform bound of norms yields
\beqn
n \in [kN,(k+1)N] \Longrightarrow \|\cL^n\psi_1-\cL^n\psi_2\|_{L^1(\fm)} \leq 2(1-a)^{n/N-1}
\eeqn
Since this holds for every $k\in\bN$, we have proved the exponential convergence 
\beq
\|\cL^n\psi_1-\cL^n\psi_2\|_{L^1(\fm)} \leq D\theta^{\alpha n}
\eeq
with parameters 
\beq
D = \frac{2}{1-a}~\textrm{ and }~\theta = (1-a)^{\frac{1}{\alpha N}}\in(0,1).
\eeq
Using parameters of the system (and the fact that $a$ must be smaller than $1/2$):
\begin{align}
D &= \frac{2}{1-1/2\exp(-(\Omega+1))} \leq 4 \\
\theta &= (1-\frac{1}{2}\exp(-(\Omega+1))^{\log(\lambda)/(4(\Omega+1))} \leq (1-\exp(-3(\Omega+1)))^{\log(\lambda)/(4(\Omega+1))}.
\end{align}

\medskip
\noindent To justify the use of the term ``coupling'' earlier, let us interpret our results in a more probabilistic fashion. First, take two initial densities $\psi^1$ and $\psi^2$. Consider random variable~$X$ and~$Y$ taking values in $\bS^1$ according to $\psi^1$ and $\psi^2$, respectively. Then, for any $n\geq 1$, $X_n = T^n(X)$ and $Y_n = T^n(Y)$ are distributed according to $\cL^n\psi^1$ and $\cL^n\psi^2$, respectively. We construct a coupling $Z_n=(X_n',Y_n')$ of~$X_n$ and~$Y_n$, in the probabilistic sense of the word, for each $n\geq 1$. To that end, notice from our construction that there exists a sequence of times $0=n_0<n_1<n_2<\cdots$ with the properties that, for any $k\geq 0$ and $n_k\leq n<n_{k+1}$,
\beq\label{eq:psi_imag}
\cL^n\psi^i = (1-a)^k\psi^i_n + (1-(1-a)^k)\rho_n,
\eeq
where $\psi^i_n$ and $\rho_n$ are probability densities, and $\rho_n$ does not depend on $\psi^i$ at all. A good definition of $Z_n$ is obtained as follows. Toss a weighted coin with the probability of heads being $1-(1-a)^k$. In the case of heads, we set $X_n' = Y_n' = R_n$ where $R_n$ has the density~$\rho_n$ on~$\bS^1$. In the case of tails, we let $X_n'$ and $Y_n'$ be distributed according to~$\psi_n^1$ and~$\psi_n^2$, respectively, independently of each other. It is easy to check that the distributions of~$X_n'$ and~$Y_n'$ are the same as those of~$X_n$ and~$Y_n$. Thus,~$Z_n$ is indeed a coupling of~$X_n$ and~$Y_n$.

\medskip
\noindent A widely known result in probability theory states that the following coupling inequality holds for \emph{an arbitrary} coupling:
\beqn
\|\cL^n\psi^1 - \cL^n\psi^2\|_{L^1(\fm)} \leq 2\bP(X_n'\neq Y_n').
\eeqn
Using the specific coupling introduced above, we have, for $n_k\leq n<n_{k+1}$, that $X_n' = Y_n' = R_n$ with probability $1-(1-a)^k$, which implies $\bP(X_n'\neq Y_n')\leq (1-a)^k$. This leads to the very same exponentially decaying estimate as in the previous section. (Of course, we could have just used~\eqref{eq:psi_imag}. The point here is that~\eqref{eq:psi_imag} was used only in the construction of the coupling~$Z_n$; had we been given~$Z_n$ directly, we could have derived the exponential bound again without~\ref{eq:psi_imag}. In any case, we have demonstrated that the work in the previous sections is genuinely related to the idea of coupling as understood by probabilists.)

\medskip
\noindent A more comprehensive, standard introduction to the coupling techniques is Torgny Lindvall's \cite{Lin}. It may also be instructive to compare our coupling argument to the one used in the classical proof of the exponential convergence of Markov chains towards a unique equilibrium distribution \cite{Ha}.

\medskip
\noindent Furthermore, we can extend our results to the general H\"older densities. This results in the desired bound for the correlation coefficients on page \pageref{consid_page}. Although this is a straighforward consequence of Theorem \ref{thm:general_density}, we establish the exact bound in Theorem \ref{thm:exp_conv} since this is the main result of the section.

\begin{theorem}[Exponential convergence of densities]~\\
\label{thm:general_density}
There exist constants $\tilde D\in\bR$ and $\theta\in(0,1)$ such that
\beqn
\|\cL^n\psi-\phi\|_{L^1(\fm)} \leq \tilde D(1+H_\alpha(\psi))\theta^{\alpha n}
\eeqn
holds for every $n\in \bN$, for every $\psi\in C^\alpha$, for every (fixed) $\alpha \in (0,1)$.
\end{theorem}

\begin{proof}
First, remark that $H_\alpha(\log(\cL^n\phi)) \leq 1+\Omega$ for sufficiently large $n$, by Corollary~\ref{cor_stric_pos} and Theorem~\ref{thm:Hoelder_estimate}. Thus, the invariance of $\phi$ results in $H_\alpha(\log(\phi)) \leq 1+\Omega \leq K$ and we can apply previous theorems to $\phi$ as well.

\medskip
\noindent We can split $\psi = 2(\psi+1)/2-1$. Thus,
\begin{align*}
\|\cL^n\psi-\phi\|_{L^1(\fm)} &= \|2\cL^n\frac{\psi+1}{2}-\cL^n\bfone-\phi\|_{L^1(\fm)} \\
&= \|2\cL^n\frac{\psi+1}{2}-2\phi-(\cL^n\bfone-\phi)\|_{L^1(\fm)}\\
&\leq 2\|\cL^n\frac{\psi+1}{2}-\phi\|_{L^1(\fm)}+\|\cL^n\bfone-\phi\|_{L^1(\fm)}\\
&\leq 2\|\cL^n\frac{\psi+1}{2}-\phi\|_{L^1(\fm)}+D\theta^{-\alpha n},
\end{align*}
where $D$ and $\theta$ are from previous results (Since $H_\alpha(\bfone) = 0$). Additionally, $(\psi+1)/2\in \cH_{H_\alpha(\psi)}$ since 
\beqn
|\log(\frac{\psi(x)+1}{2})-\log(\frac{\psi(y)+1}{2})| \leq \sup_x(\frac{2}{\psi(x)+1})|\frac{\psi(x)+1}{2}-\frac{\psi(y)+1}{2}| \leq H_\alpha(\psi) d(x,y)^\alpha.
\eeqn
That is,
\beqn
H_\alpha(\log((\psi+1)/2)) \leq H_\alpha(\psi).
\eeqn
If $H_\alpha(\psi)\leq K$ ($K$ from theorem \ref{thm:positive_convergence}), we have already proved exponential convergence of the form:
\beqn
\|\cL^n\psi-\phi\|_{L^1(\fm)} \leq 3D\theta^{\alpha n}.
\eeqn

\medskip
\noindent On the other hand, if $H_\alpha(\psi) > K$, $(\psi+1)/2\in \cH_{1+\Omega}$ after time $N(H_\alpha(\psi))=:N_0$ and by previous results
\beqn
\|\cL^n\frac{\psi+1}{2}-\phi\|_{L^1} \leq D\theta^{\alpha(n-N_0)}\qquad\forall~n>N_0.
\eeqn
If $n\leq N_0$,
\beqn
\|\cL^n\frac{\psi+1}{2}-\phi\|_{L^1} \leq \|(\psi+1)/2\|_{L^1}+\|\phi\|_{L^1} \leq 2\theta^{\alpha n}\theta^{-\alpha N_0}
\eeqn
Therefore, for every $n\in\bN$ it holds that
\beqn
\|\cL^n\frac{\psi+1}{2}-\phi\|_{L^1} \leq \max\{D\theta^{-\alpha N_0},2\theta^{-\alpha N_0}\}\theta^{\alpha n} = D\theta^{-\alpha N_0}\theta^{\alpha n},
\eeqn
if $D$ is chosen to equal $4$ as previously. Thus,
\beqn
\|\cL^n\psi-\phi\|_{L^1} \leq D(2\theta^{-\alpha N_0}+1)\theta^{\alpha n}\qquad\forall~n\in\bN.
\eeqn

\medskip
\noindent Combining the results with $\cH_{H_\alpha(\psi)} > K$ and $\cH_{H_\alpha(\psi)} \leq K$ results in
\beqn
\|\cL^n\psi-\phi\|_{L^1} \leq \max\{3,2\theta^{-\alpha N_0}+1\}D\theta^{\alpha n} = (1+2\theta^{-\alpha N_0})D\theta^{\alpha n}.
\eeqn
Plug in the appropriate values\footnote{There is no loss of generality in assuming that $H_\alpha(\log(\psi)) > 1$ since if this not the case, $\theta^{-\alpha N_0} < 1$ and the bound we obtained holds as well.}:
\begin{align*}
\|\cL^n\psi-\phi\|_{L^1} &\leq  (1+2\cdot(1-e^{-3(\Omega+1)})^{-(\log(H_\alpha(\psi))/4(\Omega+1))})4\theta^{\alpha n}
\\
& \leq 4(1 + 2\cdot 2^{(\log(H_\alpha(\psi))/4(\Omega+1))}) \theta^{\alpha n}
\\
&\leq 4(1 + 2\cdot H_\alpha(\psi)) \theta^{\alpha n}
\\
&\leq 8(1+H_\alpha(\psi))\theta^{\alpha n}.
\end{align*}
Thus, we have proved the exponential estimate for all $\psi\in C^\alpha$ with
\begin{align*}
\tilde D &= 8
\\
\theta &= (1-\exp(-3(\Omega+1))^{(\log(\lambda)/(4(\Omega+1)))}.
\end{align*}
\end{proof}

\noindent Finally, we combine the relevant results of this section and establish the exact exponential bound for the correlation coefficients.

\begin{theorem}[Exponential decay of correlation coefficients]~\\
\label{thm:exp_conv}
For bounded $f$ and H\"older continuous $g$,
\beq
\bigg|\int_X f\circ T^n g\,\rd\mu -\int_Xf\,\rd\mu\int_Xg\,\rd\mu\bigg| \leq C\|f\|_\infty(\|g\|_\infty+H_\alpha(g))\theta^{\alpha n},
\eeq
where
\beqn
C = 96(2+\Omega)^{2}
\eeqn
and $\theta$ is from the previous theorem. In particular, $\theta$ and $C$ are independent of $f$, $g$ and $\alpha$.
\end{theorem}

\begin{proof}
Considerations on page \pageref{consid_page} results in
\beqn
\bigg|\int_X f\circ T^n g\rd\mu -\int_Xf\,\rd\mu\int_Xg\,\rd\mu\bigg| \leq 3\|g\|_\infty\|f\|_\infty\left\|\cL^n\left(\frac{\phi(g+2\|g\|_\infty)}{\mu(g)+2\|g\|_\infty}\right)-\phi\right\|_{L^1(\fm)}
\eeqn
where $g \neq \bar{0}$ and $\phi$ is the invariant density. Due to Theorem \ref{existence_density}, it is evident that
\begin{align*}
&\quad H_\alpha\left(\frac{\phi(g+2\|g\|_\infty)}{\mu(g)+2\|g\|_\infty}\right) \\
&\leq \frac{H_\alpha(\phi(g+2\|g\|_\infty))}{\mu(g)+2\|g\|_\infty}\\
&\leq \frac{\|\phi\|_\infty}{\|g\|_\infty}H_\alpha(g)+3H_\alpha(\phi) \\
&\leq \frac{\mathrm{Lip}(\phi)+1}{\|g\|_\infty}H_\alpha(g)+3\mathrm{Lip}(\phi) \leq (\frac{H_\alpha(g)}{\|g\|_\infty}+3)(2+\Omega)^2
\end{align*}
and we can therefore use theorem \ref{thm:general_density} to bound the difference:
\begin{align*}
&\quad \bigg|\int_X f\circ T^n g\,\rd\mu -\int_X f \,\rd\mu \int_X g\,\rd\mu\bigg|
\\
&\leq 3\|g\|_\infty \|f\|_\infty 8 \left(1 + (\frac{H_\alpha(g)}{\|g\|_\infty} + 3)(2 + \Omega)^2)\right) \theta^{\alpha n}
\\
&\leq 96\|f\|_\infty (2 + \Omega)^{2}\left(\|g\|_\infty + H_\alpha(g)\right) \theta^{\alpha n},
\end{align*}
where $\theta$ is as earlier (see the end of Theorem \ref{thm:general_density}).
\end{proof}

\section{Mixing}

Finally, we are going to prove the mixing property of the dynamics, that is,
\beq
\label{eq:mixing2}
\int_X f\circ T^n g\rd\mu \longrightarrow \int_Xf\rd\mu\int_Xg\rd\mu\qquad\forall~f,g\in L^2
\eeq
by generalising the known result for H\"older-functions ($\alpha$ being fixed) step by step:
\begin{enumerate}
\item \eqref{eq:mixing2} holds for bounded $f$ and $g\in C^\alpha$.
\item for bounded $f$ and $g\in C^0$.
\item for bounded $f$ and simple $g$.
\item for bounded, measurable $f$ and $g$.
\item for $f,g\in L^2$.
\end{enumerate}

\noindent In what follows, asumming a function $f$ is integrable we use extensively the following abbreviated notation:
\beqn
\mu(f) = \int_Xf\,\rd\mu.
\eeqn

\noindent Step 1: We have already proved this in Theorem \ref{thm:exp_conv}.

\medskip
\noindent Step 2: $C^\alpha$ is obviously a subalgebra of $C^0$ which contains constants. Additionally, the identity map (which separates points) belongs to $C^\alpha$ and therefore, by Stone--Weierstrass theorem, $C^\alpha$ is dense in $C^0$. Thus, for every $f,g\in C^0$ and $\epsilon > 0$ there are $\tilde{f},\tilde g\in C^\alpha$ such that $\|f-\tilde f\|_\infty < \epsilon$ and $\|g-\tilde g\|_\infty < \epsilon$. Therefore
\begin{align*}
&\quad|\mu(f\circ T^n g)-\mu(f)\mu(g)| \\
&= |\mu(f\circ T^n g)-\mu(f\circ T^n \tilde g) + \mu(f\circ T^n \tilde g)-\mu(\tilde f\circ T^n\tilde g)+\mu(\tilde f\circ T^n\tilde g)\\
&\quad-\mu(f)\mu(g)+\mu(f)\mu(\tilde g)-\mu(f)\mu(\tilde g)+\mu(\tilde f)\mu(\tilde g)-\mu(\tilde f)\mu(\tilde g)| \\
&\leq \mu(|f\circ T^n||g-\tilde g|)+\mu(|f-\tilde f|\circ T^n|\tilde g|)+\mu(|f|)\mu(|g-\tilde g|)\\
&\quad +\mu(|\tilde g|)\mu(|f-\tilde f|)+|\mu(\tilde f\circ T^n\tilde g)-\mu(\tilde f)\mu(\tilde g)| \\
&\leq 2\epsilon(\|\tilde f\|_{L^1(\mu)}+\|\tilde g\|_{L^1(\mu)})+|\mu(\tilde f\circ T^n\tilde g)-\mu(\tilde f)\mu(\tilde g)| \\
&<2\epsilon(\|\tilde f\|_{L^1(\mu)}+\|\tilde g\|_{L^1(\mu)}+1)
\end{align*}
when $n$ is large enough. This proves step 2.

\medskip
\noindent Step 3: As a consequence of Lusin's theorem, for every simple $f,g$ there are continuous $\tilde f,\tilde g$ such that $\|f-\tilde{f}\|_{L^1(\fm)}<\epsilon$ and $\|g-\tilde{g}\|_{L^1(\fm)}<\epsilon$. Therefore, using the fact that the claim holds for continuous functions and $\|h\|_{L^1(\fm)} \leq \|h\|_\infty~\forall~h$ :
\begin{align*}
&\quad |\mu(f\circ T^ng)-\mu(f)\mu(g)| \\
&\leq \mu(|f\circ T^n||g-\tilde g|)+\mu(|f-\tilde f|\circ T^n|\tilde g|)+\mu(|f|)\mu(|g-\tilde g|)\\
&\quad +\mu(|\tilde g|)\mu(|f-\tilde f|)+|\mu(\tilde f\circ T^n\tilde g)-\mu(\tilde f)\mu(\tilde g)| \\
&\leq 2\epsilon(\|\tilde f\|_\infty+\|\tilde g\|_\infty)+|\mu(\tilde f\circ T^n\tilde g)-\mu(\tilde f)\mu(\tilde g)| \\
&<2\epsilon(\|\tilde f\|_\infty+\|\tilde g\|_\infty+1)
\end{align*}
when n is large enough.

\medskip
\noindent Step 4: Similar to step 2. In compact measurable spaces, bounded measurable functions can be uniformly approximated by simple functions.

\medskip
\noindent Step 5: Remark, that bounded, measurable, functions are dense in $L^2$ by simple truncation argument. Therefore, for every $f,g\in L^2$ and $\epsilon > 0$ there are bounded, measurable, $\tilde{f},\tilde g$ such that $\|\tilde f-f\|_{L^2(\fm)} < \epsilon$ and $\|\tilde f-f\|_{L^2(\fm)} < \epsilon$. Consequently,
\begin{align*}
&\quad |\mu(f\circ T^ng)-\mu(f)\mu(g)| \\
&\leq \mu(|f\circ T^n||g-\tilde g|)+\mu(|f-\tilde f|\circ T^n|\tilde g|)+\mu(|f|)\mu(|g-\tilde g|)\\
&\quad +\mu(|\tilde g|)\mu(|f-\tilde f|)+|\mu(\tilde f\circ T^n\tilde g)-\mu(\tilde f)\mu(\tilde g)| \\
&\leq 2\epsilon(\|\tilde f\|_{L^2(\fm)}+\|\tilde g\|_{L^2(\fm)})+|\mu(\tilde f\circ T^n\tilde g)-\mu(\tilde f)\mu(\tilde g)| \\
&<2\epsilon(\|\tilde f\|_{L^2(\fm)}+\|\tilde g\|_{L^2(\fm)}+1),
\end{align*}
when n is large enough. Remark, that we used Cauchy-Schwartz inequality extensively to bound $L^1$-norms.

\medskip
\noindent Finally, we wish to point out that the mixing measure we found is the only ergodic measure that has a density (w.r.t.~Lebesgue measure). First, there is a general result that ergodic measures are mutually singular \cite[p.~152]{Wal00}. Additionally, Corollary \ref{cor_stric_pos} in the Appendix states that the domain of the mixing measure covers the whole circle. By noticing that mixing measures are ergodic, together these theorems imply that if there is another ergodic measure, then its domain has a zero Lebesgue measure. That is, the measure is not absolutely continuous w.r.t.~Lebesgue measure.

\newpage
\appendix
\section{}
Two theorems, which implement very central properties of expanding dynamics, are presented here since fitting them smoothly inside the treatise turned out to be difficult. Nevertheless, that does not diminish their value.

\begin{theorem}
Suppose $\psi\in C^\alpha$ (for fixed $\alpha$). Then $\cL^n\psi\in C^\alpha$ for every $n\in\bN$ and its H\"older constant satisfies
\beq
H_\alpha(\cL^n\psi) \leq \left(\frac{H_\alpha(\psi)}{\lambda^\alpha}+(\exp(\Omega)-1)(H_\alpha(\psi)+1)\right)(1+\Omega).
\eeq
\end{theorem}
\begin{proof}
Let $|x-y|$ be small enough so that we can use the estimate from lemma \ref{lem:distortion}. Then it holds that
\begin{align*}
&\quad\;|\cL^n \psi(x)-\cL^n\psi(y)| \\
&= \bigg|\sum^{k}_{i=1}\frac{\psi(\Ti(x))}{(T^n)'(\Ti(x))}-\frac{\psi(\Ti(y))}{(T^n)'(\Ti(y))}\bigg| \\
&\leq \sum^{k}_{i=1}\frac{|\psi(\Ti(x))-\psi(\Ti(y))|}{(T^n)'(\Ti(x))}+\sum^{k}_{i=1}\bigg|\frac{\psi(\Ti(y))}{(T^n)'(\Ti(x))}-\frac{\psi(\Ti(y))}{(T^n)'(\Ti(y))}\bigg| \\
&\leq H_\alpha(\psi)\sum^{k}_{i=1}\frac{d(\Ti(x),\Ti(y))^\alpha}{(T^n)'(\Ti(x))}+\sum^{k}_{i=1}\frac{|\psi(\Ti(y)|}{(T^n)'(\Ti(y))}\bigg|\frac{(T^n)'(\Ti(y))}{(T^n)'(\Ti(x))}-1\bigg| \\
&\leq \frac{H_\alpha(\psi)}{\lambda^{\alpha n}}\sum^{k}_{i=1}\frac{d(x,y)^\alpha}{(T^n)'(\Ti(x))}+\|\psi\|_\infty\sum^{k}_{i=1}\frac{1}{(T^n)'(\Ti(y))}|e^{\Omega|x-y|}-1| \\
&\leq \frac{H_\alpha(\psi)}{\lambda^{\alpha n}}\|\cL^n\mathbf{1}\|_\infty d(x,y)^\alpha+\|\psi\|_\infty\|\cL^n\mathbf{1}\|_\infty|e^{\Omega d(x,y)}-1| \\
&\leq\left(\frac{H_\alpha(\psi)}{\lambda^\alpha}+(\exp(\Omega)-1)\|\psi\|_\infty\right)\|\cL^n\mathbf{1}\|_\infty d(x,y)^\alpha
\end{align*}
Additionally, for probability density $\psi$ there exists $x_0$ such that $\psi(x_0)=1$. Thus, for every $x$ it holds that
\beqn
|\psi(x)| = |\psi(x)-\psi(x_0)+\psi(x_0)| \leq H_\alpha(\psi)d(x,x_0)^\alpha+1 \leq H_\alpha(\psi)+1.
\eeqn
Combining this with the uniform bound for $\|\cL^n\mathbf{1}\|_\infty$ results
\beq
H_\alpha(\cL^n\psi) \leq \left(\frac{H_\alpha(\psi)}{\lambda^\alpha}+(\exp(\Omega)-1)(H_\alpha(\psi)+1)\right)(1+\Omega).
\eeq
\end{proof}

\begin{theorem}
Suppose $\psi$ is a H\"older continuous probability density. Then there exist $a>0$ and $N_1 = N_1(\psi)\in\bN$ with the property that $\inf_x \cL^n\psi(x)\geq a$ for all $n\geq N_1$.
\end{theorem}

\begin{proof}
First, remark that, by previous theorem, there is $L \in\bR_+$ such that $H_\alpha(\cL^n\psi) \leq L$ for every $n\in\bN$.

\medskip
\noindent Let's start with $\psi$. Since its probability density, there is a point $x_0\in\bS^1$ such that $\psi(x_0) = 1$. Using the fact that $H_\alpha(\psi) \leq L$ we have
\beqn
d(x,x_0) < \left(\frac{1}{2L}\right)^{1/\alpha} \Longrightarrow |\psi(x)-1| = |\psi(x)-\psi(x_0)| \leq Ld(x,x_0)^\alpha < \frac{1}{2} \Longrightarrow \psi(x) \geq \frac{1}{2}
\eeqn
More briefly:
\beqn
x \in B(x_0,(1/(2L))^{1/\alpha}) \Longrightarrow \psi(x) \geq \frac{1}{2}.
\eeqn
If $L \leq  2^{1-\alpha}$ the claim holds for every $a \in [0,1/2]$. Thus, we may assume $L > 2^{1-\alpha}$. 

Because of the expanding property of $T$, there is some $N_1$ such that
\beqn
1 < \lambda^{N_1}|B(x_0,(1/(2L))^{1/\alpha})|
\eeqn
and furthermore $T^{N_1}(B(x_0,(1/(2L))^{1/\alpha}) ) = \bS^1$. This is the case when
\beqn
\lambda^{N_1} > (2L)^{1/\alpha} > 1/|B(x_0,1/(2L)^{1/\alpha})| \textrm{ i.e. } N_1 > \frac{\log(2L)}{\alpha\log(\lambda)}.
\eeqn
Now set $N_1 =1+ \lceil \log(2L)/(\alpha\log(\lambda)) \rceil$. (We write $\lceil x \rceil$ for the smallest integer $\geq x$.) By the fact that the preimage $T^{-N_1}x$ of every point $x \in \bS^1$ intersects $B(x_0,1/(2L)^\alpha)$, we have that
\beqn
\cL^{N_1}\psi(x) = \sum_{y \in T^{-N_1}\{x\}}\frac{\psi(y)}{(T^{N_1})'(y)} \geq \frac{1}{2\|(T^{N_1})'\|_\infty} \geq \frac{1}{2\|T'\|^{N_1}_\infty}.
\eeqn
Additionally, since $H_\alpha(\cL^{n-N_1}\psi) \leq L$ for \emph{every} $n > N_1$, we can do the same reasoning and conlclude that 
\beq
\cL^n\psi(x) = \cL^{N_1}(\cL^{n-N_1}\psi)(x) \geq  \frac{1}{2\|T'\|^{N_1}_\infty}.
\eeq
This finishes the proof (by definition of infimum).
\end{proof}

\noindent The following corollary is immediate.

\begin{cor}\label{cor_stric_pos}
The invariant density $\phi$ found earlier is strictly positive.
\end{cor}



\begin{thebibliography}{2}

\bibitem{OSY} William Ott, Mikko Stenlund, and Lai-Sang Young, \emph{Memory loss for time-dependent dynamical systems}, Mathematical Research Letters, 16(3):463-475, 2009.

\bibitem{Liv} Carlangelo Liverani, Lecture notes for the course ``Sistemi Dinamici'' (2011-2012) [pdf]. Available at: $<$\url{www.mat.uniroma2.it/~liverani/SysDyn10/chap6.pdf}$>$ [Accessed 4 November 2012].

\bibitem{LSY} Lai-Sang Young, \emph{Ergodic Theory of Differentiable Dynamical Systems} [pdf], 1995. Available at: $<$\url{cims.nyu.edu/~lsy/papers/lecturenotes93.pdf}$>$ [Accessed 4 November 2012]. Also published in ``Real and Complex Dynamical Systems'', ed.~Branner and Hjorth, NATO ASI series, Kluwer Academic Publishers, (1995).

\bibitem{Lin} Torgny Lindvall, \emph{Lectures on the Coupling Method}, Wiley, New Yourk, 1992.

\bibitem{Ha} Olle H\"aggstr\"om, \emph{Finite Markov Chains and Algorithmic Applications}, Cambridge University Press, 1st edition, 2002.

\bibitem{Wal00} Peter Walters, \emph{An Introduction to Ergodic theory}, Springer, 1st softcover edition, 2000.
 
\end{thebibliography}
\end{document}